\documentclass{article}
\pdfoutput=1
\usepackage{amsmath,amssymb}
\usepackage[pdftex]{graphicx}

\newtheorem{prop}{Proposition} 
\newtheorem{theo}[prop]{Theorem}  
\newtheorem{lemm}[prop]{Lemma}    
    
\newtheorem{remark}[prop]{Remark} \let\oldrema\remark
\renewcommand{\remark}{\oldrema\normalfont}
\newtheorem{defi}[prop]{Definition}\let\olddefi\defi
\renewcommand{\defi}{\olddefi\normalfont}
\def\proof{\noindent\textit{Proof. }}
\newcounter{examp}
\def\example{\addtocounter{examp}{1}\par\bigskip\noindent\textbf{Example \arabic{examp}. }}
\newcounter{subtop}
\newcommand{\subtopic}[1]{\addtocounter{subtop}{1}\noindent\textit{\arabic{subtop}. #1. }}
 
\newcommand{\pic}[4]{\vspace{1ex}\setlength{\unitlength}{1cm}
\begin{picture}(0,#3)(5,.5)
\put(#2){\includegraphics[#4]{#1.pdf}}
\end{picture}\vspace{1ex}}

\def\C{{\mathbb C}}

\def\R{{\mathbb R}}

\def\DS{\displaystyle}
\def\qed{\ \hfill\raise3.5pt\hbox{\framebox[1ex]{\ }}}

\def\vecpopt{{\vec p}^{\,\rm opt}}
\def\vecqopt{{\vec q}^{\,\rm opt}}
\def\qst{q^*}
\def\vecqst{{\vec q}^{\,*}}
\def\pinit{{\vec p}^{\,\rm init}}
\def\qinit{{\vec q}^{\,\rm init}}
\def\Sub{\mbox{Sub}}
 
\begin{document}
%%%%%%%%%%%%%%%%%%%%%%%%%%%%%%%%%%%%%%%%%%

\vspace*{15ex}
\begin{center}\parbox{.75\textwidth}{
\begin{center}\Large
  \textbf{  Mixed Linear-Nonlinear \\Least Squares Regression
}\\[1ex] 
\large Alberto Herrera-Gomez${}^1$
\\ R. Michael Porter${}^2$

\today

 \end{center}

\noindent ${}^1$ Cinvestav-Quer\'etaro, Libramiento Norponiente 2000, Real de Juriquilla, 76230 Queretaro, Mexico, aherrerag@cinvestav.mx. Partially supported by grant CB-2012-01-179304 of  CONACyT (Mexico)\\
\noindent ${}^2$ Departamento de Matem\'aticas,
Cinvestav-Quer\'etaro,  mike@math. cinvestav.edu.mx. Partially supported by grant 166183 of CONACyT
(Mexico) }\end{center}

\bigskip
\begin{quote} \noindent\textbf{Abstract.} The problem of fitting
  experimental data to a given model function
  $f(t; p_1,p_2,\dots,p_N)$ is conventionally solved numerically by
  methods such as that of Leven\-berg-Marquardt, which are based on
  approximating the $\chi^2$ measure of discrepancy by a quadratic
  function.  Such nonlinear iterative methods are usually necessary
  unless the function $f$ to be fitted is itself a linear function of
  the parameters $p_n$, in which case an elementary linear
  Least Squares regression is immediately available. When linearity is
  present in some, but not all, of the parameters, we show how to
  streamline the optimization method by reducing the ``nonlinear
  activity'' to the nonlinear parameters only.  Numerical examples are
  given to demonstrate the effectiveness of this approach.  The main
  idea is to replace entries corresponding to the linear terms in the
  numerical difference quotients with an optimal value easily obtained
  by linear regression.  More generally, the idea applies to
  minimization problems which are quadratic in some of the parameters.
  We show that the covariance matrix of $\chi^2$ remains the same even
  though the derivatives are calculated in a different way. For this
  reason, the standard non-linear optimization methods can be fully
  applied.
\end{quote}



\section{Introduction}

Linear least squares regression provides a global, fast, and absolute
minimum of Chi-squared when the function to be fitted depends linearly on
all its parameters. However, when the dependence on even a single
parameter fails to be linear, linear regression can no longer be
applied, and the solution must be successively approximated. The
classical approach, which is widely employed, is to use non-linear
least squares regression, which is a procedure treating all the
variables as non-linear. In this paper we show how to take advantage
of linear regression when the dependence on some of the parameters is
linear.  The benefits of the method here described are many; in
particular it is more robust and frequently faster than the classical
approach. In addition, all the benefits of the non-linear methods
employing gradients and their covariance matrix still apply.

The techniques which have been developed since the time of I.\ Newton
for calculating local minima of a smooth function
$F(p_1,p_2,\dots,p_N)$ of several real variables form a fundamental
part of numerical analysis and have been refined to improve
performance as much as possible.  They derive from the well known fact
that the solution of the minimization problem for a quadratic function
of the parameters $p_1,p_2,\dots,p_N$, say,
\[ \sum_{n,n'} A_{n,n'}p_np_{n'} + \sum_n b_np_n + C, \] reduces to
finding vector quantity $(1/2)A^{-1}\vec b$, i.e., to simply solving a
system of linear equations, to determine the point at which the
gradient of the function vanishes. In real life, $F$ is usually not of such
a particularly simple form, and one must turn to nonlinear methods
\cite{Ac,DB,PTVF}; these are commonly based on approximating $F$ by a
 a quadratic function near a supposed minimum $\vec p$, and improving
it via estimations of the gradient $\vec b$ and the Hessian
matrix $A$ of $F$ at $\vec p$. 

We are particularly interested in the case when the function $F$ to be
optimized (e.g., $\chi^2$ in the case of fitting experimental data
with a model $f$) is expressible as a quadratic polynomial in some,
but not all, of the variables (for example if we add a single
transcendental term to the formula above). Then the optimization
problem must be considered as essentially nonlinear.  The purpose of
this article is to show how to take advantage of as much of the
quadratic structure of $F$ as may be present, in order to
significantly improve the optimization procedure in many cases.

An important application of optimization principles is to the theory
of fitting formulas to data.  Given a model function
$f_{p_1,p_2,\dots,p_N}(t)$ of a variable $t$, and a set of
experimental data values $(t_1,y_1),(t_2,y_2),\dots$, one seeks the
combination $\vec p$ of the parameters $p_n$ for which
$f_{\vec p}(t_j)$ is most nearly equal to $y_j$.  The $\chi^2$
measure of the discrepancy from an exact fit, seen as a function of
the $p_n$, indeed becomes a quadratic polynomial when $f_{\vec p}(t)$
is linear in all these parameters (see section \ref{sec:fit}).  Linear
models occur in many situations, such as polynomial or Fourier
approximations, and have the pleasant characteristic of reducing to
least squares problems which can be solved by very fast algorithms
\cite{Bev}.

As an example, an approximation widely used, both in experimental
sciences and in theoretical mathematical studies, is the exponential
fitting
\[ f_{p_1,p_2,\dots,p_{3N}}(t) = 
  \sum_n p_{3n}  e^{-\left(\frac{p_{3n-1}}{t-p_{3n-2}}\right)^2 }. 
\] 
This is linear in those $p_n$ for which the index $n$ is a multiple
of 3.  By reducing the number of parameters treated nonlinearly by one
third, not only do we obtain a reduction in computation time in many
cases, but the numerical stability of the fitting problem can also be
greatly improved. The same holds for models with alternatives to the
classical exponential peaks.   
 
The method given here was first applied to fit infrared \cite{HVM} and
photoemission \cite{HHM} spectra, and was implemented in the software
\textit{AAnalyzer} \cite{Her-AA1} in 1998.  More information about
this software can be found elsewhere \cite{Her-AA2}.  Although the
fundamental idea of our method is extremely simple, it has not been
analyzed mathematically until now; to our knowledge,  
all other commercially available software treats fitting problems
as totally nonlinear when so much as one nonlinear term is present.

In Section \ref{sec:shortcut} we introduce the notion of ``shortcut
derivative'' of a differentiable function and explain its use for
accelerating optimization methods.  In Section \ref{sec:fit} we apply
this concept for acceleration of virtually any standard fitting
algorithm.  A theoretical result justifying the of the use of this
technique is proved in Section \ref{sec:correlation}, and some
numerical examples are given in Section \ref{sec:numerical}. In the
final discussion section we compare the efficiency of our method with
standard methods.

\section{Shortcut derivatives\label{sec:shortcut}}

We present here the concept of ``shortcut derivative'' in a fairly
general context, to be specialized later, and explain its relevance
to optimization problems.

\subsection{Optimization relative to a subset of  parameters}
Consider a smooth real-valued function $F(\vec p,\ \vec q)$ to be
optimized, where we have arbitrarily separated the variables into a
length-$M$ vector $\vec p=(p_1,p_2,\dots,p_M)$ is of real numbers
varying in some region $\Omega_1$ of Euclidean space $\R^M$ and an
$N$-vector $\vec q=(q_1,q_2,\dots,q_N)\in\Omega_2\subseteq\R^N$. Let
$1\le m\le M$.  For $\delta>0$, we modify the classical difference
quotient
\begin{eqnarray}  \label{eq:numderiv}
 D_F(\delta,m) &=& \frac{1}{2\delta}\big(F(\vec p + \delta e_m,\vec q)\  
   -    F(\vec p - \delta e_m,\vec q) \big),
\end{eqnarray}
which is the standard numerical approximation for the partial
derivative $\partial F(\vec p,\vec q)/\partial p_m$, as follows.  Here
we have written $\vec e_m=(0,0,\dots0,1,0,\dots0)$ for the $m$-th
canonical basis vector of $\R^M$, i.e.,
$\vec p + \delta e_m= (p_1,p_2,\dots,p_{m-1},p_m+\delta,p_{m+1}\dots,p_M)$.

\defi \label{def:shortcut} For each fixed $\vec p$ define
$\vecqst(\vec p)\in\R^N$ as the value of $\vec q$ minimizing
$F(\vec p,\vec q)$, that is, satisfying
\begin{equation}  \label{eq:qstdef}
    F(\vec p,\vecqst(\vec p)) = \min_{\vec q\in\Omega_2}  f (\vec p,\vec q)   .
  \end{equation}
Then the \emph{shortcut derivative} of
  $F(\vec p,\ \vec q)$ with respect to the single parameter $p_m$, and
  relative to the parameter subset $\vec q$, is the limit of
\begin{eqnarray}  \label{eq:shortcutderiv}
  D^*_F(\delta,m)
  &=& \frac{1}{2\delta}\big(F(\vec p + \delta e_m,\,\vecqst(\vec p + \delta e_m))
      - F(\vec p - \delta e_m,\, \vecqst(\vec p - \delta e_m))\big)  
\end{eqnarray}
as $\delta\to0$. Similarly we have the \emph{second shortcut
  derivative} as the limit of
 \begin{eqnarray}  
 D^*_F(\delta,m,m')
  &=&   \frac{1}{4\delta^2} \bigg( F( p+\delta\vec e_m+\delta\vec e_{m'},\,
           \vecqst( p+\delta\vec e_m+\delta\vec e_{m'}))  \nonumber\\
  && \quad\quad  - F(\vec p-\delta\vec e_m+\delta\vec e_{m'},\,
           \vecqst( p-\delta\vec e_m+\delta\vec e_{m'}))  \nonumber\\[1ex]
  && \quad\quad   - F(p+\delta\vec e_m-\delta\vec e_{m'},\,
           \vecqst( p+\delta\vec e_m-\delta\vec e_{m'}))  \nonumber\\
  && \quad\quad   + F(p-\delta\vec e_m-\delta\vec e_{m'},\,
           \vecqst( p-\delta\vec e_m-\delta\vec e_{m'}))  \bigg) \label{eq:shortcutderiv2}
\end{eqnarray}
 
When discussing shortcut derivatives, we will assume, as is common in
studies of numerical techniques, that the minimizer $\vecqst(\vec p)$
  is unique. This may often be achieved by
working locally, i.e., by reducing $\Omega_1$ to a region of
interest. The existence of the limits (\ref{eq:shortcutderiv}),
(\ref{eq:shortcutderiv2}), will be established in the following discussion.

The gradient derivative of $F$ separates naturally into two parts,
\begin{equation}  \label{eq:gradpgradq}
  \nabla F = (\nabla_pF,\ \nabla_qF),
\end{equation}
where $\nabla_p=(\partial/\partial p_1,\dots,\partial/\partial p_M)$
and $\nabla_q=(\partial/\partial q_1,\dots,\partial/\partial
q_N)$. Thus for fixed $\vec p$, the vector $\vecqst(\vec p)$ can be
characterized by the property
\begin{equation}  \label{eq:qstarcondition} 
  \nabla_qF|_{\vecqst(\vec p)} = 0.
\end{equation}
It follows from the Implicit Function Theorem that
$\vecqst\colon\Omega_1\to\R^N$ is then a smooth function, assuming
that the Jacobian matrix of the correspondence
$\vec q\mapsto \nabla_qF|_{(\vec p,\,\vec q)}$ is nonsingular.
 
The notion of shortcut derivative is intimately connected with the
``reduced function'' $F^*\colon\Omega_1\to\R$, which we define by
\begin{equation}  \label{eq:F*}
  F^*(\vec p) = F(\vec p,\,\vecqst(\vec p) ).
\end{equation}
Applying the Chain Rule for derivatives to (\ref{eq:F*}) yields
\begin{equation}  \label{eq:partialsF*}
  \frac{\partial F^*}{\partial p_m} =
  \frac{\partial F}{\partial p_m} + \sum_{n=1}^{N}
    \frac{\partial F}{\partial q_n} \frac{\partial q^*_n}{\partial p_m},
\end{equation}
where the right-hand side is evaluated at $(\vec p,\vecqst(\vec p))$,
and  $q^*_n\colon\Omega_1\to\R$ are the coordinate functions of $\vecqst$,
\[ \vecqst(\vec p) = (q^*_1(\vec p),\dots,q^*_N(\vec p)). \] As an
immediate consequence of (\ref{eq:qstarcondition}) and
(\ref{eq:partialsF*}), we have
\begin{prop}
  $\nabla F^*|_{\vec p} =   \nabla_p F|_{(\vec p,\vecqst(\vec p))}$.
\end{prop}
 
Further, by placing $F^*$ in place of $F$ in (\ref{eq:numderiv}) one
immediately sees the following.

\begin{prop} \label{pr:main} The shortcut derivative of $F(\vec p,\vec
  q)$ with respect to $p_m$, relative to the variables of $\vec q$, is equal to
  $\partial F^*(\vec p)/\partial p_m$.
  \end{prop}
 
\begin{remark}\label{rem:Fstar}
  Obviously, when $(\vecpopt,\vecqopt)$ is a minimum point for $F$,
  necessarily $\vecqopt=\vecqst(\vecpopt)$, and $\vecpopt$ is a
  minimum point for $F^*$. Conversely, a minimum point $\vecpopt$ for
  $F^*$ generates a minimum point $(\vecpopt,\vecqopt)$ for $F$ via
  $\vecqst$. This fact is relevant under the assumption that
  \emph{minimum points of $\vecqst$ may be obtained at low
    computational cost}.  In general, this will hold when the first
  and second partial derivatives of $F$ with respect to the variables
  $q_n$ may be calculated at low computational cost.  In particular
  this holds when $F$ is a quadratic polynomial in $q_1,\dots,q_n$ for
  every fixed $\vec p$.
\end{remark}

\subsection{Shortcut acceleration of optimization methods}
With the above ingredients we can already outline the  
the shortcut algorithm. It is simply the minimization of $F^*$ of equation
(\ref{eq:F*}).

Specifically, consider any of the well-known methods of optimization
which require calculation or estimation of the gradient vector
$\nabla$ and the Hessian matrix $H$ of the function $F$ to be
optimized. Assume that it is feasible to calculate the function
$\vecqst$ defined implicitly by (\ref{eq:qstdef}).  Begin with an
initial guess $\pinit$, for the parameters $p_1,\dots,p_N$, and
evaluate
\begin{equation}  \label{eq:calcq*}
   \vecqst(\pinit\pm\delta_m\vec e_m) 
\end{equation}
for $1\le m\le M$. The offsets $\delta_m$ may be chosen according to
the sensitivity of $F$ in each variable.  The values (\ref{eq:calcq*})
are applied in (\ref{eq:shortcutderiv}) and (\ref{eq:shortcutderiv2})
to approximate the shortcut partial derivatives of $F$ with respect to
$p_m$, which are used in place of the usual numerical approximations
of the derivatives to form the shortcut gradient $\nabla^*_p$ and the
shortcut Hessian matrix $H^*$, which are then used in place of the
true gradient and Hessian in the chosen algorithm---this will be
justified in Section \ref{sec:correlation}. The algorithm produces an
improved value for the minimizer $\vec p$, which then as usual takes
the place of $\pinit$ in the following iteration.  When sufficient
accuracy is obtained, a final application of $\vecqst$ to the
resulting $\vec p$ completes the desired optimal parameter $(\vec p,
\vecqst(\vec p))$ as explained in Remark \ref{rem:Fstar}.

\section{Fitting of mixed linear-nonlinear models\label{sec:fit}}

We apply the above considerations to explain how to manage the problem
of fitting experimental data to a model function which is linear in
some variables $q_1,\dots,q_N$ and (possibly) nonlinear in the
remaining variables $p_1,\dots,p_M$. Such a function can be expressed
in the general form
\begin{equation}  \label{eq:modelfunction}
   f(t)= f_{\vec p,\vec q}(t) = \sum_{n=1}^N q_n \varphi_{n,\vec p}(t) + \psi_{\vec p}(t). 
\end{equation}
for some functions $\varphi_{n,\vec p}$ ($n=1,\dots,N$) and
$\psi_{\vec p}$ which do not depend on $\vec q$. (The term
$\psi_{\vec p}$ rarely appears in physical applications, and the
reader may wish to ignore it in what follows.) One wishes to choose
$\vec p,\vec q$ so as to minimize (perhaps locally) the Chi-squared
quantity \cite{Ac,DB,PTVF}
\begin{equation}  \label{eq:chisq}
 \chi^2 = \chi^2(\vec p,\vec q) = \sum_t  (f_{\vec p,\vec q}(t)-y_t)^2 w_t
\end{equation}
This sum is taken over a finite collection of sample values of the
independent variable $t=t_1,t_2,\dots$, to which there are associated
measurements $y_t$, and   respective weights $w_t>0$. It is common
to take $w_t=y_t^{-1}$ as the inverse of the covariance of $y_t$, but the
particular choice of $w_t$ will not be relevant to our considerations.

We now consider the minimization of $F=\chi^2$ in the context of the
previous section. The first and second partial derivatives of
(\ref{eq:chisq}) with respect to $p_m$ are
\begin{equation}  \label{eq:chisqdp1}
\frac{\partial(\chi^2y)}{\partial p_m} =
      -2\sum_t (f(t) - y_t) \frac{\partial f(t)}{\partial p_m} w_t,
\end{equation}
\begin{equation}  \label{eq:chisqdp2}
 \frac{\partial^2(\chi^2)}{\partial p_m\partial p_{m'}} =
    2\sum_t  \left(
      \frac{\partial f(t)}{\partial p_m}
      \frac{\partial f(t)}{\partial p_{m'}}
      - (f(t) - y_t)
       \frac{\partial^2 f(t)}{\partial p_m\partial p_{m'}}
 \right)  w_t .
\end{equation}
As discussed in \cite{PTVF}, the second terms of the summands in
(\ref{eq:chisqdp2}) are generally discarded in numerical work not only
because the factors $f(t) - y_t$ tend to be small, but because
greater numerical stability is achieved this way. Thus for fitting problems
we will make no use of (\ref{eq:shortcutderiv2}).

The  form  (\ref{eq:modelfunction}) for $f$ leads to
the expansion
\begin{eqnarray}  
  \chi^2(\vec p,\vec q) &=& \nonumber
  \sum_n\sum_{n'}\left(\sum_t \varphi_n\varphi_{n'} w_t\right)q_nq_{n'}  
  + \sum_n\left(\sum_t2 (\psi-y_t)\varphi_nw_t\right)q_n\\ 
 && \ \ + \sum_t (\psi-y_t)^2w_t , \label{eq:chisqquadratic}
\end{eqnarray}
which is quadratic in $q_1,\dots,q_N$ for fixed $\vec p$.
With this, the partial derivatives of $\chi^2$ with respect to the
linear parameters $q_n$ are easily obtained in terms of
\begin{equation}  \label{eq:fdq1}
  \frac{\partial f(t)}{\partial q_n} \ = \ \varphi_n(\vec p,t) .
\end{equation}
Namely, one finds that (\ref{eq:chisqdp1}) becomes
\begin{equation} \label{eq:chisqdq1}  
  \frac{\partial(\chi^2)}{\partial q_n} =  2\sum_{n'}\left(\sum_t  \varphi_n\varphi_{n'}w_t\right)q_{n'} + 2\sum_t  (\psi-y_t)\varphi_nw_t
\end{equation}
which conveniently represents the gradient of $\chi^2$ in the form
$\nabla\chi^2 = A\vec q + \vec b$. Thus for fixed $\vec p$, the
(absolute) minimum of $\chi^2$ is attained when $A\vec q + \vec
b=0$. Typically there are more data points $t$ than the number $N$ of
parameters and this is an overdetermined linear system, so there does
not exist an exact solution to this system, but the residual
$A\vec q + \vec b$ is minimized in terms of the $L_2$ norm by the
linear least squares regression \cite{PTVF} which is available as a
standard function of many numerical software packages.  In this sense,
by Definition \ref{def:shortcut}, this linear least squares solution is the
best approximation for $\vecqst(\vec p)$.  In contrast to the
$p_m$-derivatives, the quantities involved in setting up this linear
system for the $q_n$-derivatives require no special numerical
derivation, since they are already calculated whenever
$f_{\vec p,\vec q}(t)$ itself is evaluated.

\section{Reduced covariance matrix\label{sec:correlation}}
 
The Hessian matrix $H$ appearing in a nonlinear optimization procedure
is associated with the correlations of the fitted model with the
original data; we will study it here and comment more fully in Section
\ref{sect:comparison} below.

We return to the generality of a function $F_{\vec p,\vec q}$ to be
minimized as in Section \ref{sec:shortcut}.  Let $H^*(\delta)$ denote
the approximation of the shortcut Hessian matrix corresponding to a
parameter displacement $\delta>0$. (More precisely, $\delta=\delta_m$
refers to a displacment in a single direction $\vec e_m$ for
notational simplicity; our statemens will be valid as well for a
vector of displacements $(\delta_1\dots\delta_M)$.)  The true Hessian
matrix $H$ of $F$ admits a natural block decomposition
\[   H = \left(  \begin{array}{cc}
     H_{pp} &   H_{pq} \\   H_{qp} &   H_{qq} 
  \end{array}  \right) \in \R^{M+N,M+N} 
\]
with $H_{pp}\in\R^{M,M}$, $H_{pq}\in\R^{M,N}$, $H_{qp} \in \R^{N,M}$,
and $H_{qq} \in \R^{N,N}$. The symmetries $H_{pp}=( H_{pp})^{\rm T}$,
$H_{qp}=(H_{pq})^{\rm T}$, $H_{qq}=(H_{q})^{\rm T}$ are a consequence
of the fact that the order of differentiation is irrelevant. Consider
the inverse matrices
\begin{eqnarray*}
   \boldsymbol{\eta}(\delta) &=& \left( \eta^*_{jj'}(\delta) \right) =    (H^*(\delta))^{-1}, \\
   \boldsymbol{\eta} &=& (H_{pp})^{-1},
\end{eqnarray*}
and denote
\[  H^*(0) = \lim_{\delta\to0} H^*(\delta) . \]

It is well known that the matrix entries of $H^{-1}$ represent the
covariances of the full set of fitted parameters
$p_1,\dots,p_M,q_1, \dots,q_N$, and the diagonal elements, the
variances, are of particular importance (see for example \cite[Sect.\
15.4, 15.6]{PTVF}, especially equation (15.4.15)). Similarly, the
entries of $\eta$ are the covariances of the parameters
$p_1,\dots,p_M$ for the fitting corresponding to the minimization of
(\ref{eq:F*}).  The following result was discovered empirically.

\begin{theo} \label{th:Hessian} Suppose that the Hessian matrix $H$ is
  invertible. Then $H^*({\delta})$ is also invertible for each
  $\delta>0$, and
 \[   \eta^*_{m,m'}(\delta) \to \eta_{m,m'} \mbox{ as } \delta\to0
\] for $1\le m\le M$, $1\le m'\le M$. In other words,
 $H^*(\delta)^{-1}\to H_{pp}^{-1}|_{(p,q^*)}$ as $\delta\to0$.
\end{theo}

We devote the rest of this section to the proof of this result, which
depends on two lemmas. Note that we do not claim that $H^*(\delta)$
approximates $H$ in any way.

\begin{lemm}\label{lemm:1}
  The Jacobian matrix of the function $\vecqst$,
\begin{equation}  \label{eq:Jq*}
   J\vecqst = \left( \frac{\partial \qst_n}{\partial p_m}
     \right)_{\substack{\!\!1\le n\le N\\\!\!\!1\le m\le M}} \in \R^{N\times M} 
   \end{equation}
   is
   given by the formula
\begin{equation}
  J\vecqst = - H_{qq}^{-1}H_{qp} .
\end{equation}
\end{lemm}

\proof From (\ref{eq:qstarcondition}) and the Chain Rule,
\[  0 = \frac{\partial^2 F}{\partial q_n\partial p_m}
   (\vec p,\vecqst(\vec p))  = \left.
   \frac{\partial F}{\partial q_n\partial p_m} + \sum_{n'=1}^N
          \frac{\partial^2 F}{\partial q_n\partial q_{n'}} 
          \frac{\partial\qst_{n'}}{\partial p_m} 
 \ \right\vert_{\vec q=\vecqst(\vec p)}  
\] 
for $1\le m\le M$, $1\le n\le N$. In matrix form this is
$0=H_{qp}+H_{qq}(J\vecqst)$. The invertibility of $H$ implies the 
invertibility of $H_{qq}$, from which the result follows.  \qed

\begin{lemm}\label{lemm:2} $H^*(0) = H_{pp}- H_{pq}H_{qq}^{-1}H_{qp}$.  
\end{lemm}

\proof
Again by the Chain Rule,
\[ \frac{\partial F^*}{\partial p_m} (\vec p) = 
     \frac{\partial F}{\partial p_m}(\vec p,\vecqst(\vec p)) + \sum_{n=1}^N
          \frac{\partial F}{\partial q_n}(\vec p,\vecqst(\vec p)) 
          \frac{\partial\qst_n}{\partial p_m}(\vec p,\vecqst(\vec p)) 
\]
for $1\le m\le M$. Now differentiate with respect to $p_{m'}$,
\begin{eqnarray*}
 \frac{\partial^2 F(\vec p,\vecqst(\vec p))}{\partial p_m\partial p_{m'}}   &=& 
 \left(  \frac{\partial F }{\partial p_m\partial p_{m'}} + \sum_{n=1}^N
          \frac{\partial^2 F }{\partial p_m\partial q_n}%\cdot 
          \frac{\partial \qst_n }{\partial p_{m'}} 
 \ \right)  \\  && \hspace*{-8ex}
   + \sum_{n=1}^N \left[ \left(\frac{\partial F }{\partial p_{m'}\partial q_n}
       + \sum_{n'=1}^N
          \frac{\partial^2 F }{\partial q_n\partial q_{n'}}
          \frac{\partial \qst_{n'} }{\partial p_{m'}} 
 \ \right)\frac{\partial \qst_n}{\partial p_{m}} + 
 \frac{\partial F}{\partial q_n} 
  \frac{\partial^2\qst_n}{\partial p_m\partial p_{m'}}
 \right] \\ 
 &=&   \ \frac{\partial F }{\partial p_m\partial p_{m'}} 
    + \sum_{n=1}^N  \left(
           \frac{\partial^2 F }{\partial p_{m}\partial q_n} 
          \frac{\partial \qst_n }{\partial p_{m'}}  +
        \frac{\partial^2 F }{\partial p_{m'}\partial q_n} 
          \frac{\partial \qst_n }{\partial p_{m}}  \right) \\
  && + \   \sum_{n=1}^N\sum_{n'=1}^N  \frac{\partial^2 F }{\partial q_n\partial q_{n'}} 
  \frac{\partial \qst_{i}}{\partial p_{m}}\frac{\partial \qst_{n'}}{\partial p_{m'}} 
    \ + \  \sum_{n=1}^N  \frac{\partial F}{\partial q_n} 
              \frac{\partial^2\qst_n}{\partial p_m\partial p_{m'}}     
\end{eqnarray*} 
for $1\le m'\le M$. By (\ref{eq:qstarcondition}), the
last sum vanishes, and upon substituting the definitions of the entries of
$H_{pp}$, $H_{pq}$, $H_{qq}$, we find
 \[   h^*_{mn}(0) = h_{m,m'} 
    + \sum_{n=1}^N  \left(
       h_{mm'}\frac{\partial \qst_n }{\partial p_{m'}}  +
       \frac{\partial \qst_n }{\partial p_m} h_{nm'}  \right) 
   + \   \sum_{n=1}^N\sum_{n'=1}^N  
   \frac{\partial\qst_n}{\partial p_m} h_{nn'}
   \frac{\partial\qst_{n'}}{\partial p_{m'}} .
 \]
 In matrix notation this is
 \[  H^*(0)  =    H_{pp} + H_{pq}J_{\vecqst}+ (J_{\vecqst})^{\rm T}H_{qp}
  +  (J_{\vecqst})^{\rm T}H_{qq} J_{\vecqst},
\]
so upon substituting  the formula of Lemma \ref{lemm:1} and canceling, we arrive
at the result. \qed

\noindent\textbf{Proof of Theorem \ref{th:Hessian}. }
By Cramer's rule, we want to compare the matrix entries 
\[ \begin{array}{rcccl}
  \eta_{m,m'}      &=& \DS  \frac{\det \Sub_{m,m'}H}{\det H}, \\[2ex]
  \eta^*_{m,m'}(0) &=& \DS  \frac{\det \Sub_{m,m'}H^*(0)}{\det H^*(0)}.
\end{array}
\]
where we write $\Sub_{m,m'}A$ for the submatrix of $A$ obtained by removing
the $m$-th row and $m'$-th column. It will suffice to show
\begin{eqnarray}
    \det H^*(0) &=& \frac{1}{\det H_{qq}}\det H, \label{eq:Hmat}
  \\     \det \Sub_{m,m'}H^*(0) &=& \frac{1}{\det H_{qq}}\det \Sub_{m,m'}H, \label{eq:Hsubmat}
\end{eqnarray}
since upon dividing these two formulas we obtain $\eta_{m,m'}^*(0)=\eta_{m,m'}$.

We deduce (\ref{eq:Hmat})  from Lemma \ref{lemm:2} and the
formula for the determinant of a block matrix
\[   \det  \left(  \begin{array}{cc}
     A & B\\C&D   \end{array} \right) =
    \det D\det(A-BD^{-1}C),    
\]
valid when $A$ and $D$ are square submatrices and  $D$ is invertible. Indeed,
\[  \det H  = \det H_{qq}\ \det(H_{pp}- H_{pq}H_{qq}^{-1}H_{qp}) 
  =  \det H_{qq} \ \det H^*(0).
\]
To verify (\ref{eq:Hsubmat}) we apply $\Sub_{m,m'}$ to both sides of the
formula of Lemma \ref{lemm:2},
\[  \Sub_{m,m'}H^*(0) = \Sub_{m,m'}H_{pp}- \Sub_{m,m'}[H_{pq}H_{qq}^{-1}H_{qp}].
\]
Applying the general rule $\Sub_{m,m'}(ABC)=(\Sub_{m0}A)B(\Sub_{0m'}C)$, we find the
block structure
\[  \Sub_{m,m'} H = \left(  \begin{array}{cc}
    \Sub_{m,m'}H_{pp}  & \Sub_{m0}H_{pq}\\ \Sub_{0m'}H_{qp} &   H_{qq} \end{array} \right) 
\]
and invoking Lemma \ref{lemm:2} again, the block determinant is
\begin{eqnarray*}
 \det  \Sub_{m,m'} H &=& \det  H_{qq} 
     \det\bigg( \Sub_{m,m'}H_{pp} - (\Sub_{m0}H_{pq})(H_{qq}^{-1})(\Sub_{0m'}H_{qp}) \bigg)  \\ &=& \det  H_{qq} 
     \det\bigg( \Sub_{m,m'}H_{pp} -  \Sub_{m,m'}[H_{pq}H_{qq}^{-1}H_{qp}] \bigg)  \\[1ex]
   &=& \det  H_{qq} \det  \Sub_{m,m'} H^*(0),
\end{eqnarray*}
which is (\ref{eq:Hsubmat}). This completes the proof. \qed
 
 \section{Numerical examples\label{sec:numerical}}
 
We compare the results of solving fitting problems in a model which
is linear in the variables $q_1,\dots,q_N$ by   the
method of shortcut derivatives and by the classical approach, which
treats all variables equally, i.e., using a parameter set
\[ p_1,\dots,p_M,p_{M+1},\dots,p_{M+N} \]
where we denote $p_{M+n}=q_n$ for $1\le n\le N$.  The conversion to
the classical formulation is easily programmable in terms of the
expression (\ref{eq:modelfunction}) simply by incorporating the
$\varphi_n$ terms into the function $\psi$ and reindexing the
variables.

\begin{figure}[b!]\centering %1a
\pic{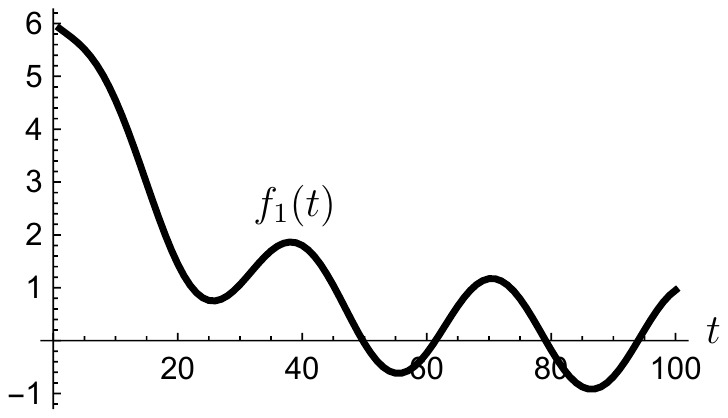}{-2.3,-.2}{4.4}{scale=.9} 
\pic{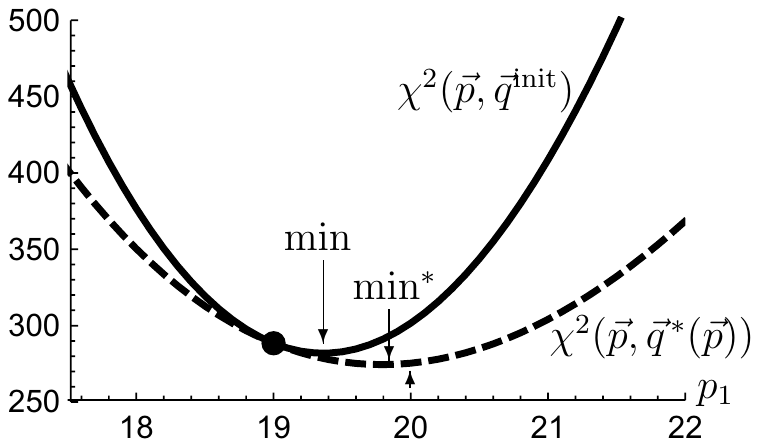}{4.2,-1}{4.}{scale=1} 
\caption{\label{fig:better}\small Model curve (\ref{eq:2p2q}) for
  ``true'' values of its 4 parameters (left).  With $p_2$ fixed, plots
  of $\chi^2$ are given (right) as a function of $p_1$ where the
  initial optimum value $\qinit$ is fixed (solid), and with $\vec q$
  optimized for each $p_1$ (dashed). The minima of these plots are
  approximations for the true minimum at $p_1=20$.}
\end{figure}
 
In the following numerical experiments, the method of
Levenberg-Marquardt \cite{Mar,PTVF} was programmed in
\textit{Mathematica} (Wolfram) and run on an ordinary laptop
computer. We avoided fancy variations such as in \cite{TS}.

\example %1
Consider a model function with four parameters, 
 \begin{equation}\label{eq:2p2q}
  f_1(t) =  q_1e^{-t/p_1} +q _2 \sin\frac{t}{p_2}, \quad t=1,2,\dots,100, 
\end{equation}
with the particular values $(p_1,p_2)=(20,5)$, $(q_1,q_2)=(6,1)$ as in
Figure \ref{fig:better}.

 \begin{figure}[t!]\centering %1b
\pic{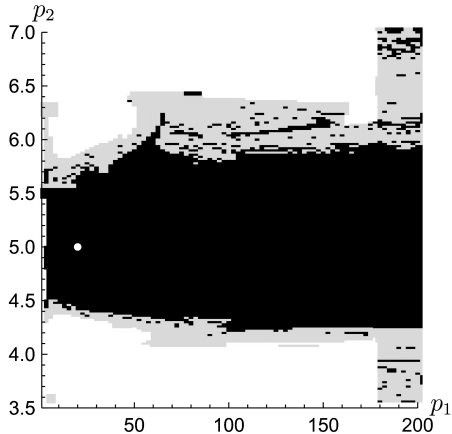}{1.7,-.9}{5}{scale=.8} 
\caption{\label{fig:sensitivity}\small Comparison of shortcut method
  to classical fitting for Example 1. In the black region, both
  methods converge to the true value (marked with a white spot). The
  outer white area is the region of $p_1$,$p_2$ values for which
  neither method converges to the correct $f_1(t)$, and in the intermediate
  gray area only the shortcut method works. The many small
  black ``islands'' in the gray area  contain values where the
  classical method may be considered to converge ``accidentally'';
  one could not reliably choose an initial guess near such points. }
\end{figure}

The nearby parameter $\pinit=(19,4.9)$ is chosen for illustrative
purposes.  The parameter set is automatically completed with
$\qinit=\vecqst(p^{\rm init})=(6.19664, 0.947731)$ approximately. We
consider the question of minimizing $\chi^2$ as $p_1$ varies while
leaving $p_2=4.9$ fixed, to see how close we can come back to the true
value pf $p_1=20$. There are two natural ways to do this:
(a) to leave the value of $\vec q$ fixed, thus considering parameter sets
$(p_1, p_2^{\rm init}, q_1^{\rm init}, q_2^{\rm init})$ in which only $p_1$
varies; (b) to optimize the $q$-values along with $p_1$, thus considering
$(p_1, p_2^{\rm init}, \vecqst(p_1, p_2^{\rm init}))$. Note that case (a) is
simply a straight line in the 4-dimensional parameter space, while
(b) is a curve in this space specially adapted to our problem.

One sees in Figure \ref{fig:better} that the graph of the values of
$\chi^2$ for case (a) lies above the graph of $\chi^2$ for case
(b). The minimum value of the former curve is approximately 19.35,
whereas the minimum of the latter is approximately 19.85, much closer
to the true value $p_1=20$. This phenomenon illustrates why the
shortcut method tends to require fewer iterations to approximate
minima; of course, in general one is locating a minimum on a
higher-dimensional hypersurface of which we have been able to illustrate
only a 1-dimensional slice here.
  
For the same model function (\ref{eq:2p2q}), Figure
\ref{fig:sensitivity} depicts the region in the $p_1$-$p_2$ plane of
initial guesses from which both the shortcut method and the classical
method converge to the correct values $p_1=20$, $p_2=5$, lying
properly within the region where only the shortcut method gives the
correct solution.  In order to make this comparison, the classical
method was started with the ``guess values'' of $q_1,q_2$ given by
$\vecqst(p_1,p_2)$.  In practice, one would almost surely make a poorer
initial guess if $\vecqst$ were not evaluated, and the contrast
between the sensitivity of the two methods to the initial guess would
be found to be even greater than shown here.

\example %2
For this example we take as model function a variable number of Gaussian peaks, 
\begin{equation}  \label{eq:almostlinearmodel}  
  f_2(t) =   e^{-(t/p_1)^2}  + \sum_{n=1}^N q_n e^{-(t-n)/5)^2} .
\end{equation}
The peaks centered at the positive integer points $t=1,2,\dots,N$ have fixed 
widths but variable heights.  In contrast, the single peak centered at $t=0$
has variable width determined by $p_1$, the only nonlinear parameter, sufficient
to preclude the sole use of least squares regression. 

 \begin{figure}[h!]\centering
\pic{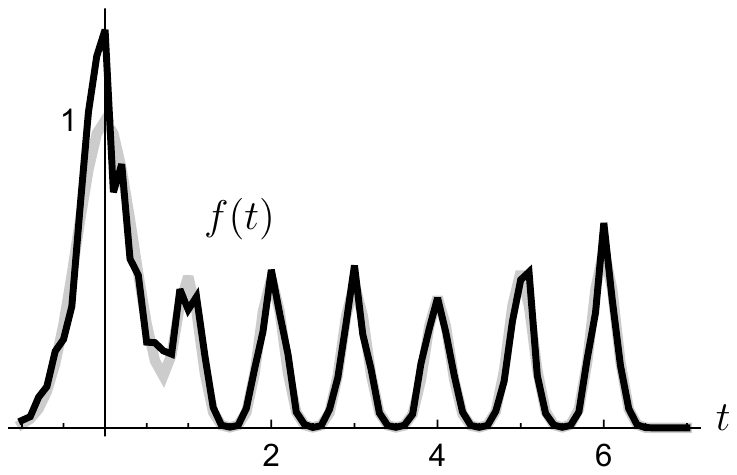}{-2,-1}{4.}{scale=.8}\pic{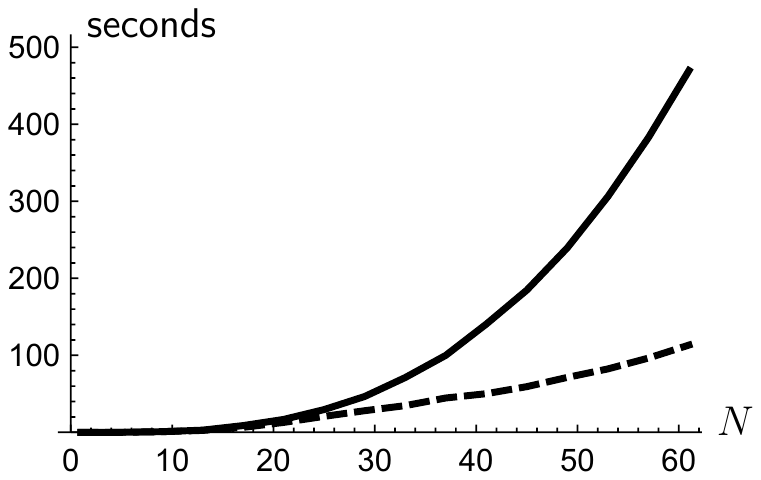}{4,-1.3}{4}{scale=1}
\caption{\label{fig:almost}\small Sample of model function
  (\ref{eq:almostlinearmodel}) (left)
  in gray, with randomized data points ($\pm30\%$) superimposed in black.  Running
  times (right) averaged over 5 runs for each case.}
\end{figure}

Figure \ref{fig:almost} shows the advantage of the shortcut method for
this example. The increasing number of data points (with a spacing
$\Delta t=0.1$ between consecutive points) is an additional factor in
requiring more calculation time as the number $N+1$ of peaks
increases.  It is seen that the running times are approximated by
$0.01\times N^{2.3}$ and $0.003\times N^{2.9}$. For the shortcut
method, the Levenberg-Marquardt required 2 iterations for the entire
range of values of $N$, while the classical method required 3
iterations for $N>45$. The number of function calls in the classical
method grew steadily from $50\%$ greater for low values of $N$ to 35 times
greater for $N=60$, easily offsetting the higher cost of the calculation
of $\vecqst$ in each derivative.

\example %3

A situation which presents a greater opportunity for the shortcut
method is the simultaneous fitting of several ``files'' of data.  An
application occurs in X-ray photoelectron spectroscopy (XPS)
\cite{MH}, in which the chemical composition of a sample material is
to be determined by the energy distribution of the electrons leaving a
surface illuminated by X-rays at different angles (one file for each
angle). For physical reasons the centers and widths of the peaks are
not affected by a change in angle of the incident beam, so these
parameters are common  to the collection of model functions which must be fit
simultaneously (``shared-parameters hypothesis'', as described
in \cite{MH}).

To simulate this phenomenon we use the model function
\begin{equation}  \label{eq:filesmodel}
  f_3(t) = q_1e^{((t-p_1)/p_4)^2} +  q_2e^{((t-p_2)/p_5)^2} +  q_3e^{((t-p_3)/p_6)^2} + q_4t+q_5 
\end{equation} 
in which the last two terms represent the ``background noise''.  (In
practice many other models are often used in which the peaks are not
Gaussian.)  Let us suppose that sample data is given for
$k=1,2,\dots,K$ readings, sharing common values of $p_1,\dots,p_6$ but
each with its own set of linear parameters
\[  \vec q^{(k)} = (q_{k,1},\dots,q_{k,5})
\]

\begin{figure}[t!]\centering  %3a
  \pic{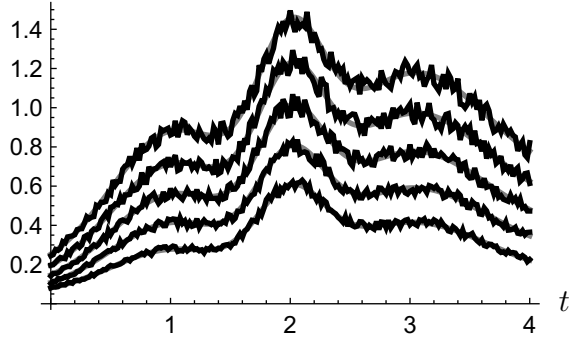}{1,-1.2}{4.}{scale=1.}
  \caption{\label{fig:5files}\small Simulation of five files of data results
    formed of three slightly overlapping Gaussian peaks, with increasing
    heights approximately proportional to the noise baseline. Here $0<t<4$.
  }
\end{figure}

\begin{figure}[b!]\centering %3b
\pic{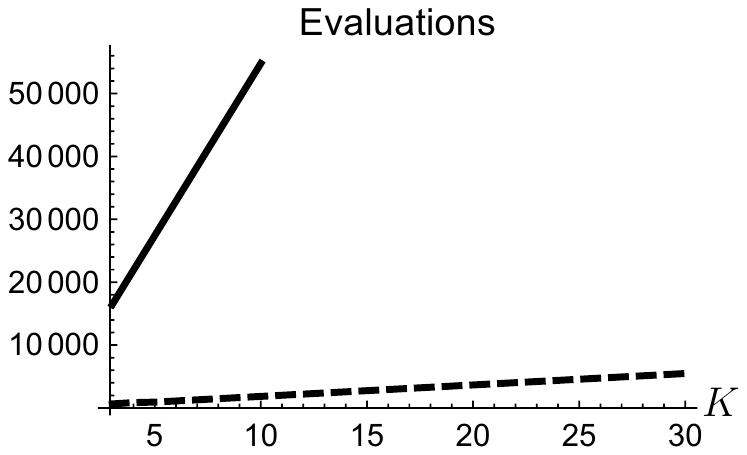}{-2,-1.3}{4}{scale=.95}
\pic{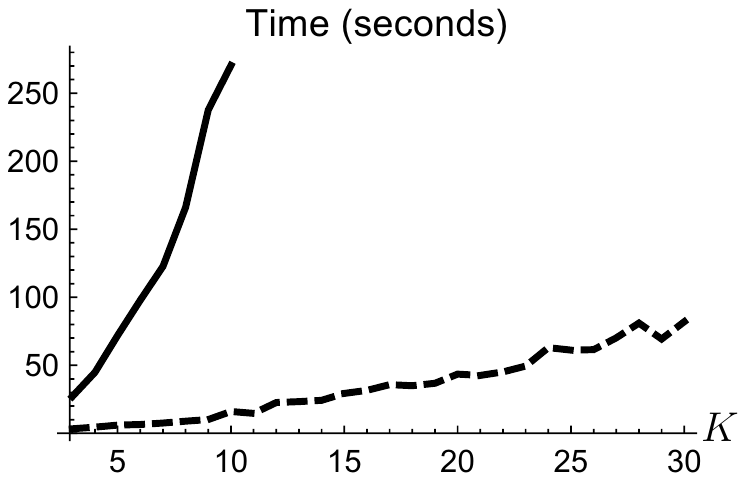}{4.8,-1.3}{4}{scale=.95}
\caption{\label{fig:filesresults}\small Count of function evaluations
  and computing time for shortcut algorithm (dashed) and traditional
  Levenberg-Marquardt calculation (solid) for increasing numbers $K$ of
  simultaneous files of data modeled by (\ref{eq:filesmodel}).   }
\end{figure}

As illustrated in Figure \ref{fig:5files}, there is a set of sample
data $\vec y^{(k)}$ for each $k$, which determines in turn a
discrepancy function $\chi^2_k(\vec p,\vec q)$. The closeness of the
collective fitting corresponding to the $5K+6$ parameters $\vec p$,
$\vec q^{(1)}$, \dots, $\vec q^{(K)}$ is measured by the sum
\begin{equation} \label{eq:Phi}
  \Phi(\vec p,\vec q^{(1)},\vec q^{(2)},\dots,\vec q^{(K)}) =
  \sum_{k=1}^K \chi^2_k(\vec p,\vec q).  
\end{equation}  
It should be noted that the partial derivative of $\Phi$ with respect
to $p_m$ is the sum of the corresponding derivatives of $\chi^2_k$,
while the partial derivative of $\Phi$ with respect to $q_n$ is
given by the partial derivative of the single summand of (\ref{eq:Phi}) in
which the parameter $q_n$ appears.  For this reason it is a straightforward
matter to represent the gradient and Hessian of $\Phi$ in terms of
the entries of the gradients and Hessians of the various $\chi^2_k$.
 
Using this information, we calculated the best fit via the
Levenberg-Marquardt method for the function
$\Phi^*(\vec p)= \Phi(\vec p,\vecqst{}^{(1)}(\vec
p),\vecqst{}^{(2)}(\vec p),\dots,\vecqst{}^{(2)}(\vec p))$, where
$\vecqst{}^{(k)}$ denotes the optimization of the five $q$-values
\[ \vec q^{(k)} =(q^{(k)}_1,q^{(k)}_2,\dots,q^{(k)}_5) \]
with respect to the data $\vec y^{(k)}$.  For comparison, the
classical Levenberg-Marquardt was applied after translating the
$t$-variable of the model functions in such a way as to form a single
function on an interval $K$ times as long as the original one, and
incorporating all the variables in a single $p$-list via the relation
$\vec q^{(k)} =(p_{5k-4}, p_{5k-3},\dots,p_{5k})$. The $\chi^2$ of
this auxiliary function is precisely the sum in (\ref{eq:Phi}). Thus the
resulting minima are identical.  As Figure \ref{fig:filesresults}
shows, the computational cost is reduced significantly by the use of
shortcut derivatives.
  
 \section{Comparison of methods\label{sect:comparison}}

We discuss briefly some of the differences in the results produced
by the shortcut method with respect to traditional regression methods.
  
\subsection{Direct advantages of the linear-nonlinear method}

\subtopic{Robustness of computation} As illustrated by Example 1 and
many other examples we have calculated, the linear-nonlinear method
tends to offer more flexibility in the choice of initial guess for the
optimization procedure.

The robustness of is more strongly manifested when shared parameters
are involved such as in Example 3. This phenomenon was previously
studied in the case of files of data in \cite[Fig.\ 7]{MH}, where it
was noted that precision of the assessment of the shared-parameters
($p$-variables) increases with the number of files.  Although this is,
in fact, expected because each file adds information, it could be
somewhat counterintuitive. A belief holds in parts of the physics
community that there is a limit on the number of parameters, holding
covariances among each other, that can be simultaneously optimized;
limits of 17 to 20 have been suggested.  In this example, we have
shown that it is possible to optimize a far greater number of
parameters simultaneously (e.g., 156 parameters when $N=30$) even
though the covariance between each pair of variables is not zero.  In
general, the use of the mixed linear-nonlinear method described here
tends to increase the robustness of the optimization process since
linear regression provides the best possible value of the linear
parameters at each step.  The method inherits some important
characteristics of linear regression since it is also possible to
catch the culprit parameters if the minimization has multiple
solutions (e.g., too many free parameters or too many peaks).

\subtopic{Estimation of uncertainty}
As is well known, $H_{pp}^{-1}$ evaluated at optimal $(p,q)$ is a very
relevant matrix inasmuch as its diagonal elements measure the
covariance of the $p$-parameters [13]. In the discussion of the
particular case study of \cite[Fig.\ 7]{MH} it was found numerically that
the standard deviation predicted through $H^{-1}$ is quite similar to
the actual standard deviation obtained through many trials.  Thus
Theorem \ref{th:Hessian} allows one to use the diagonal elements of
${H^*(\delta)}^{-1}$ to report the uncertainty on the $p$-parameters
even though the derivatives are obtained in a nonstandard way.

The advantage of this is that in the covariance analysis of a fitting
problem, the inversion of the smaller matrix $H^*(\delta)$ for fixed,
small $\delta>0$ is less costly than the inversion of the full matrix
$H$.

\subtopic{Operation count and computation time}
It must be recognized that each iteration of the shortcut method is
computationally more costly than the corresponding classical method,
due to the evaluations of $\vecqst$ in each shortcut derivative.

Consider the computing time required for the classical method, i.e.,
for a sequence of model functions of the form $f(p_1,\dots,p_M)$ with
increasing numbers $M$ of the parameters, the evaluation of which we
will assume implies a computational cost proportional to $M$. (This
does not hold for Example 2, where the number of data points $t$ also
grows with $M$.)  The cost of evaluating first and second partial
derivatives numerically grows as $O(M)$, which implies that the
gradient vector of $f$ costs $O(M^2)$. 
Forming the Hessian of $\chi^2$ via (13) requires $O(M^2)$ operations
once the first partial derivatives are known (since we are assuming
the number of $t$ values is fixed).  The solution of an $M\times M$
matrix equation by least squares costs $O(M^3)$, which is thus the
cost of a single iteration of an $M$-parameter fitting. 

Now we consider a model function of the form (\ref{eq:modelfunction}).
Each evaluation of $\varphi_1(\vec p),\dots,\varphi_1(\vec p)$ costs
$O(MN)$ since each function individually costs $O(M)$.  An evaluation
of $f(\vec p,\vec q)$ thus also costs $O(MN)$ (this is independent of
whether we are speaking of a single $t$ or the entire fixed set of $t$
values).  To evaluate all the partial derivatives
$\partial(\chi^2)/\partial q_n$, $1\le n\le N$ by
(\ref{eq:modelfunction}) requires $O(N^2)$ multiplications, so the
cost of the Hessian of $\chi^2$ with respect to $\vec q$ costs
whichever is greater of $O(MN)$ and $O(N^2)$. Following Definition
\ref{def:shortcut}, we see that an evaluation of $\vecqst(\vec p)$ is
carried out by minimizing with respect to the $N$ variables
$q_1,\dots,q_N$, which means solving an $N\times N$ system.  Hence the
cost of $\vecqst(\vec p)$ is $O(MN)+O(N^3)$.

The shortcut derivatives $\partial f/\partial p_m$, $1\le m\le M$ are
obtained from (\ref{eq:shortcutderiv}) with $f$ in place of $F$, with
a total cost of $O(M^2N)+O(MN^3)$, which by (\ref{eq:chisqdp2}) is the
cost of the shortcut Hessian of $\chi^2$, an $M\times M$ matrix.
Since the solution of the $N\times N$ system is only $O(N^3)$, the
iteration cost for the shortcut method is $O(M^2N)+O(MN^3)$. A more
refined analysis, for which the number $T$ of data points is allowed
to grow, gives the value $O(M^2NT)+O(MN^3)$.

In comparison, the classical method has a cost of $O((M+N)^3)$. Thus
it can be seen that the shortcut method will be more costly in the
long run when $M$ is at least of the same proportion as $N$.  In
Examples 2 and 3, the ratio $M/N$ effectively tends to zero, which
accounts for the considerable savings of time with the shortcut method.

 In our application of the Levenberg-Marquardt procedure, in many cases
we have found sometimes, but not always, that fewer iterations are
necessary than with the classical method, but never more.  In
experiments with model functions such as
\[ \sum_n q_n e^{-\left(\frac{p_{2n-1}}{t-p_{2n}}\right)^2 }, \] in
which the $p$, $q$ variables are in proportion of 2:1, with the peaks
centered at integer points $p_{2n}=1,2,\dots$,  the
shortcut method greatly reduces the number of iterations required, but
the computation time grows faster than in the traditional method, due
to the cost of calculating $\vecqst$ in the shortcut derivatives.
It would  be interesting to look for ways of reducing the cost
of $\vecqst$ if possible.

\subsection{Other approaches to linear-nonlinear regression problems}

We have reviewed the mathematical literature in optimization and
regression fairly carefully, as well as descriptions of available
software, and to the best of our knowledge, the mixed linear-nonlinear
approach does not seem to have been considered previously apart from
its use in the software \cite{Her-AA1}, where shortcut derivatives
were applied without theoretical justification.  Also, Theorem
\ref{th:Hessian} was discovered empirically \cite{MH} by noticing that
the first $M$ diagonal elements of $H^{-1}$ appeared to be numerically
equal (within rounding errors) to the diagonal elements of
$H^*(\delta)^{-1}$. The algorithm in \cite{Her-AA1} calculates the
computationally cheaper quantities $H^*(\delta)^{-1}$ as
approximations of the covariances of the parameters. The excellent
results obtained compared to other similar programs, and the
unexplained coincidence of the inverse diagonal elements motivated the
present investigation of the mathematical properties of this approach.

In fact, there has been little systematic research into the idea of
combining linear and nonlinear aspects of fitting problems. 
An ad-hoc method for $f(t)=q_1e^{p_1 t} + q_2e^{p_2t}$, suitable for
working out by hand calculation, is described in \cite{Han}.
Exponential regressions of an arbitrary number summands are discussed
in \cite{SCh}.
Chapter 9 of \cite{DW} gives a detailed discussion of a
linear-nonlinear problem, one of the few we have found on this
subject:  a process of heat produced by cement hardening 
which is nonlinear in time and linear in some of the other control
variables. The approach there alternates linear and nonlinear approximation,
but is quite different from the method described here.  Fitting of
parametrized curves $(y(t),z(t))$ in the plane is discussed in \cite{Sp},
in which linearity also plays an important role.

Many types of industrial problems (see for example \cite{SMKKS}), as
well as calculations in mathematical biology \cite{MCh,Rab},
mathematical finance, and other areas, require fitting of model curves
or surfaces to observed data.  We believe that many such areas could
benefit from the shortcut method of optimization.

\end{document}